# Entropy production by implicit Runge-Kutta schemes

## Carlos Lozano

### Computational Aerodynamics Group

### National Institute for Aerospace Technology (INTA), Spain

**Abstract**: This paper follows up on the author's recent paper "Entropy Production by Explicit Runge-Kutta schemes" [1], where a formula for the production of entropy by fully discrete schemes with explicit Runge-Kutta time integrators was presented. In this paper, the focus is on implicit Runge-Kutta schemes, for which the fully discrete numerical entropy evolution scheme is derived and tested.

## 1. Introduction

The purpose of this paper is to extend the work of [1], which computed the entropy produced by fully discrete numerical schemes with explicit Runge-Kutta (RK) integrators, to general RK (including implicit) schemes (see [2, 3] for an introduction to RK schemes with further references). The fully discrete schemes that will be analyzed must be understood as deriving from space and time discretization of general systems of conservation laws via the method of lines. The spatial discretization must be supplemented with a suitable time integrator, usually a Runge-Kutta scheme (see [4] for a recent review). While higher-order time integrators are usually explicit RK methods, implicit and semi-implicit RK schemes are also gaining relevance as viable candidates for the integration of time-dependent pde (see for example [5, 6] and the recent work by Jameson [7] and references therein). Despite their increased numerical cost, implicit RK methods offer superior stability properties, higher-order accuracy and relaxed CFL restrictions compared with explicit methods, and are ideally suited to the integration of stiff problems.

Many physical systems of practical interest are equipped with at least one entropy function that is exactly conserved in smooth solutions but is dissipated across singularities such as shocks. At the numerical level, schemes that reproduce this behavior are called entropy stable and constitute an active area of research dating back to the works of Lax [8], Harten et al. [9] and Osher [10]. Entropy stability has proved instrumental in the analysis of numerical approximations of systems of conservation laws, particularly in several space dimensions, as it provides global stability estimates for numerical methods for multi-dimensional conservation laws (see [11] for a recent review with an extended set of references). The systematic construction of entropy stable schemes for systems of conservation laws was initiated by Tadmor [12] by blending entropy conservative fluxes with appropriate dissipation operators in order to achieve entropy stability. Following this path, high-order entropy conservative fluxes were developed in [13], while arbitrarily high-order entropy stable schemes were first developed in [14]. The stability analysis of Discontinuous Galerkin schemes was initially addressed in [15] using entropy considerations (see [16, 17, 18, 19] for updated accounts and further references). Entropy stability considerations can also be used to guide the design of schemes based on summation-by-parts operators (see [20] and references therein).

In recent times, entropy-stable schemes have risen to prominence [21, 22, 23] especially in the context of turbulent flow simulations using high-order methods [24]. While significant attention is usually paid to the entropy production properties of the spatial schemes, comparatively little effort has been devoted to the study of the entropy production of temporal schemes, which however can have a significant impact [25, 26]. Indeed, the spatial entropy production can be made negative with a suitable entropy stable scheme, but time integration schemes may still produce entropy. To attain fully discrete local entropy stability (the requirement that entropy should be dissipated within a time-step at each node of the computational mesh), the spatial dissipation must dominate the temporal entropy production, which may require a suitable CFL condition (see e.g. [12]). Low-order, explicit RK schemes tend to produce (spurious) entropy, while higher-order schemes tend to dissipate it [27, 1]. The situation with implicit schemes is however less explored. Tadmor [12] (see also [28]) has shown that the Backward-Euler scheme dissipates entropy, while a generalized form of the Crank-Nicolson scheme introduced in [13] is entropy conservative[1]. For those cases, entropy stability of the fully discrete scheme is attained regardless of time-step size as long as the spatial scheme is entropy dissipative. In the general case, it will be shown that even implicit RK schemes generally produce entropy and, thus, entropy stability can only be guaranteed under a CFL bound.

The paper is organized as follows. Section 2 contains the heart of the developments of the paper. Section 2.1 derives the fully discrete entropy evolution equation for generic Runge-Kutta time integrators coupled to finite-volume spatial discretizations of systems of conservation laws in 1D. Sections 2.2 and 2.3 explore the consequences of the entropy evolution equation for local and global entropy stability, respectively, while section 2.4 discusses several RK schemes in detail. Section 3 presents the results of numerical experiments involving entropy stable discretizations of the linear advection and the Burgers equation (with quadratic and logarithmic entropies) coupled to several RK schemes. For illustration purposes, the investigation focuses in particular on the two and three stage Gauss and Radau IIA schemes, as well as on two Singly-diagonal (SDIRK) schemes analyzed in [7]. Finally, Section 4 presents the conclusions.

**2. Entropy production of fully discrete schemes**

In this paper we will primarily consider systems of conservation laws in 1D of the form

$$\frac{\partial U}{\partial t} + \frac{\partial f(U)}{\partial x} = 0 \tag{1}$$

---

[1] In its original form, this generalized Crank-Nicolson scheme requires an intermediate temporal state that unfortunately does not generally have a closed form and requires quadrature. An explicit construction has been recently derived in [25].

where $U(x,t) = (U_1(x,t),\ldots,U_N(x,t))^T$ is the vector of conservative variables and $f(U) = (f_1(U),\ldots,f_N(U))^T$ is the flux vector. We assume that the above system is endowed with an entropy pair $(\eta(U), \phi(U))$, where the entropy function $\eta$ is convex and the entropy flux $\phi$ satisfies the compatibility condition $\eta_U(U) f_U(U) = \phi_U(U)$ (the subscript $U$ denotes differentiation with respect to $U$). The entropy pair satisfies the following conservation equation [29, 30]

$$\frac{\partial \eta}{\partial t} + \frac{\partial \phi(U)}{\partial x} \leq 0 \tag{2}$$

where equality holds for smooth fields $U$. The left-hand side (LHS) of (2), which is non-vanishing for discontinuous (weak) solutions (e.g., shocks), is the entropy production.

## 2.1. Entropy evolution equation

At the discrete level, entropy also evolves with time in a way that resembles (2). The precise numerical entropy evolution equation follows from the fully discretized version of (1), that we will assume to be approximated by a semi-discrete Finite Volume numerical scheme on a grid over the range $i = 1,\ldots,n$ as

$$\Delta x_i \frac{dU_i}{dt} + F_{i+\frac{1}{2}} - F_{i-\frac{1}{2}} = 0 \tag{3}$$

where $\Delta x_i$ is the length of cell $i$ and $F_{i\pm\frac{1}{2}}$ are the numerical fluxes. Multiplying the left-hand side (LHS) of (3) by the vector of entropy variables $v_i^T = \eta_U(U_i)$ and rearranging yields the semi-discrete numerical entropy scheme

$$\Delta x_i \frac{d\eta_i}{dt} + \Phi_{i+\frac{1}{2}} - \Phi_{i-\frac{1}{2}} = \Pi_i \tag{4}$$

where

$$\Phi_{i\pm\frac{1}{2}} = \overline{v}_{i\pm\frac{1}{2}}^T F_{i\pm\frac{1}{2}} - \overline{\Theta}_{i\pm\frac{1}{2}} \tag{5}$$

is the numerical entropy flux and

$$\Pi_i = \frac{1}{2}\left(\Pi_{i+\frac{1}{2}} + \Pi_{i-\frac{1}{2}}\right)$$
$$\Pi_{i\pm\frac{1}{2}} = \Delta_{i\pm\frac{1}{2}} v^T F_{i\pm\frac{1}{2}} - \Delta_{i\pm\frac{1}{2}} \Theta \tag{6}$$

is the numerical entropy production. $\Theta = v^T f - \phi$ is the entropy potential verifying $f(U(v)) = \nabla_v \Theta$, while $\overline{\Theta}_{i+\frac{1}{2}} = \frac{1}{2}(\Theta_{i+1} + \Theta_i)$ and $\Delta_{i+\frac{1}{2}}\Theta = \Theta_{i+1} - \Theta_i$ denote averaging and differencing across the face, respectively. The right-hand side (RHS) of (4) is the numerical entropy production of the spatial scheme, $S_i^{(x)} = \Pi_i$.

For numerical computation, Eq. (4) is coupled to a time advancement scheme which will also produce (or dissipate) entropy in an amount $S^{(t)}$ that will depend on the type of scheme [12]. The previous paper [1] addressed this problem for explicit Runge-Kutta schemes, while now we will address the issue of general (i.e. implicit) Runge-Kutta schemes with Butcher's Tableau [2]

$$\begin{array}{c|ccc} c_1 & a_{11} & \cdots & a_{1s} \\ \vdots & \vdots & \ddots & \vdots \\ c_s & a_{s1} & \cdots & a_{ss} \\ \hline & b_1 & \cdots & b_s \end{array}$$

which corresponds to the following time-advancement scheme

$$U_i^{n+1} - U_i^n = -b_1 R_i^{(1)} - \cdots - b_s R_i^{(s)} \tag{7}$$

where $R_i^{(k)} = \lambda_i (F_{i+\frac{1}{2}}^k - F_{i-\frac{1}{2}}^k)$ are the numerical residuals, $F_{i\pm\frac{1}{2}}^k = F_{i\pm\frac{1}{2}}(U^{(k)})$ are the numerical fluxes, $\lambda_i = \Delta t_i / \Delta x_i$, where $\Delta t_i$ is the (possibly) spatially varying time-step, and the intermediate states $U^{(k)}$ are computed as

$$\begin{aligned} U_i^{(1)} - U_i^n &= -a_{11} R_i^{(1)} - \cdots - a_{1s} R_i^{(s)} \\ &\vdots \\ U_i^{(s)} - U_i^n &= -a_{s1} R_i^{(1)} - \cdots - a_{ss} R_i^{(s)} \end{aligned} \tag{8}$$

In order to obtain an expression for the entropy evolution of the above scheme, we multiply both sides of (7) on the left by $(v_i^{n+1})^T = v(U_i^{n+1})^T = \eta_U(U_i^{n+1})$. Using (4), the RHS can be cast as

$$\begin{aligned} -(v_i^{n+1})^T \sum_{k=1}^{s} b_k R_i^{(k)} &= -\sum_{k=1}^{s} b_k (v_i^{(k)})^T R_i^{(k)} - \sum_{k=1}^{s} b_k (v_i^{n+1} - v_i^{(k)})^T R_i^{(k)} = \\ -\lambda_i \sum_{k=1}^{s} b_k (v_i^{(k)})^T (F_{i+\frac{1}{2}}^k - F_{i-\frac{1}{2}}^k) &- \sum_{k=1}^{s} b_k (v_i^{n+1} - v_i^{(k)})^T R_i^{(k)} = \\ -\lambda_i \sum_{k=1}^{s} b_k \left( \Phi_{i+\frac{1}{2}}^k - \Phi_{i-\frac{1}{2}}^k \right) + \frac{1}{2} \lambda_i \sum_{k=1}^{s} b_k \left( \Pi_{i+\frac{1}{2}}^k + \Pi_{i-\frac{1}{2}}^k \right) &- \sum_{k=1}^{s} b_k (v_i^{n+1} - v_i^{(k)})^T R_i^{(k)} \end{aligned} \tag{9}$$

where $v_i^{(k)} = v(U_i^{(k)})$. For the LHS, the following identity due to Tadmor [12]

$$(v_i^{n+1})^T (U_i^{n+1} - U_i^n) = \eta_i^{n+1} - \eta_i^n + B_i^{n+\frac{1}{2}} \tag{10}$$

is used, where

$$B_i^{n+\frac{1}{2}} = \int_{-1/2}^{1/2} \left(\tfrac{1}{2}-\xi\right)(\Delta^{n+\frac{1}{2}}v_i)^T H(v(\xi))\Delta^{n+\frac{1}{2}}v_i \, d\xi$$

$$H(v) = \frac{\partial U}{\partial v} = (\eta_{UU})^{-1}$$

$$v(\xi) = \overline{v}_i^{n+\frac{1}{2}} + \xi \Delta^{n+\frac{1}{2}} v_i \qquad (11)$$

$$\overline{v}_i^{n+\frac{1}{2}} = \tfrac{1}{2}(v_i^{n+1} + v_i^n)$$

$$\Delta^{n+\frac{1}{2}} v_i = v_i^{n+1} - v_i^n$$

Combining (9) and (10) yields

$$(v_i^{n+1})^T (U_i^{n+1} - U_i^n) = \eta_i^{n+1} - \eta_i^n + B_i^{n+\frac{1}{2}} = -(v_i^{n+1})^T \sum_{k=1}^{s} b_k R_i^{(k)} =$$
$$-\lambda_i \sum_{k=1}^{s} b_k \left(\Phi_{i+\frac{1}{2}}^k - \Phi_{i-\frac{1}{2}}^k\right) + \frac{1}{2}\lambda_i \sum_{k=1}^{s} b_k \left(\Pi_{i+\frac{1}{2}}^k + \Pi_{i-\frac{1}{2}}^k\right) - \sum_{k=1}^{s} b_k (v_i^{n+1} - v_i^{(k)})^T R_i^{(k)} \qquad (12)$$

which can be arranged as

$$\eta_i^{n+1} - \eta_i^n + \lambda_i \sum_{k=1}^{s} b_k \left(\Phi_{i+\frac{1}{2}}^k - \Phi_{i-\frac{1}{2}}^k\right) =$$
$$-B_i^{n+\frac{1}{2}} - \sum_{k=1}^{s} b_k (v_i^{n+1} - v_i^{(k)})^T R_i^{(k)} + \frac{1}{2}\lambda_i \sum_{k=1}^{s} b_k \left(\Pi_{i+\frac{1}{2}}^k + \Pi_{i-\frac{1}{2}}^k\right) \qquad (13)$$

If the matrix $A_{ij} = a_{ij}$ is invertible[2], Eq. (8) can be used to write $R_i^{(k)}$ in terms of $U_i^{(j)} - U_i^n$ to cast (13) as

$$\eta_i^{n+1} - \eta_i^n + \lambda_i \sum_{k=1}^{s} b_k \left(\Phi_{i+\frac{1}{2}}^k - \Phi_{i-\frac{1}{2}}^k\right) =$$
$$\underbrace{-B_i^{n+\frac{1}{2}} + \sum_{k=1}^{s}\sum_{j=1}^{s} b_k A_{kj}^{-1}(v_i^{n+1} - v_i^{(k)})^T (U_i^{(j)} - U_i^n)}_{S^{(t)}} + \underbrace{\frac{1}{2}\lambda_i \sum_{k=1}^{s} b_k \left(\Pi_{i+\frac{1}{2}}^k + \Pi_{i-\frac{1}{2}}^k\right)}_{S^{(x)}} \qquad (14)$$

This is our main result, from which the entropy produced by the fully discrete scheme (7) can be computed. The bottom row of (14) is the local entropy production of the fully discrete scheme (denoted as $S_i$), which can be separated into spatial $S_i^{(x)}$ and temporal $S_i^{(t)}$ parts.

**The Backward-Euler scheme**

---

[2]If $A$ is not invertible we can consider the enlarged matrix $\tilde{A} = \begin{pmatrix} A \\ b^T \end{pmatrix}$ (where $b = (b_1,...,b_s)^T$). If $\text{rank}(\tilde{A}) = s$, we can invert (7) in the least-squares sense as $R_i = -(\tilde{A}^T \tilde{A})^{-1} \tilde{A}^T \Delta U_i$, (where $\Delta U_i = (U_i^{(1)} - U_i^n,...,U_i^{(s)} - U_i^n)^T$ and $R_i = (R_i^1,\cdots,R_i^s)^T$), so in (14) we should replace $A^{-1}$ with $(\tilde{A}^T \tilde{A})^{-1} \tilde{A}^T$.

We can check (14) with a known example such as the Backward-Euler (BE) scheme

$$U_i^{n+1} - U_i^n = -\lambda_i (F_{i+\frac{1}{2}}^{n+1} - F_{i-\frac{1}{2}}^{n+1}) \tag{15}$$

which corresponds to (7) with $s = 1$ and $a_{11} = b_1 = 1$. The entropy production of this scheme has been computed by Tadmor [12] as

$$\eta_i^{n+1} - \eta_i^n + \lambda_i (\Phi_{i+\frac{1}{2}}^{n+1} - \Phi_{i-\frac{1}{2}}^{n+1}) = -B_i^{n+\frac{1}{2}} + \lambda_i \Pi_i^{n+1} \tag{16}$$

(where Tadmor's notation has been adapted to ours) which does coincide with (14) for $s = 1$ and $a_{11} = b_1 = 1$. The RHS of (16) is the local entropy production of the BE scheme. The temporal entropy production $S_i^{(t)} = -B_i^{n+\frac{1}{2}}$ is semidefinite-negative, which implies that, for entropy stable spatial schemes[3] (for which $\Pi_i \leq 0$), the Backward-Euler scheme is entropy dissipative.

In the general case, the sign of the temporal entropy production is harder to analyze, and we will find that, in general, implicit RK schemes may produce entropy. In order to proceed further, it is advisable to cast (14) in an alternative form more amenable to subsequent analysis. Multiplying both sides of (7) on the left by $(v_i^n)^T$ and proceeding as above, only that this time a different Tadmor's identity is used on the LHS,

$$(v_i^n)^T (U_i^{n+1} - U_i^n) = \eta_i^{n+1} - \eta_i^n - E_i^{n+\frac{1}{2}} \tag{17}$$

where

$$E_i^{n+\frac{1}{2}} = \int_{-1/2}^{1/2} (\xi + \tfrac{1}{2})(\Delta^{n+\frac{1}{2}} v_i)^T H(v(\xi)) \Delta^{n+\frac{1}{2}} v_i d\xi \tag{18}$$

yields

$$\begin{aligned} S_i &= \eta_i^{n+1} - \eta_i^n + \lambda_i \sum_{k=1}^{s} b_k \left( \Phi_{i+\frac{1}{2}}^k - \Phi_{i-\frac{1}{2}}^k \right) = \\ & E_i^{n+\frac{1}{2}} - \sum_{k=1}^{s} \sum_{j=1}^{s} b_k A_{kj}^{-1} (v_i^{(k)} - v_i^n)^T (U_i^{(j)} - U_i^n) + \frac{1}{2} \lambda_i \sum_{k=1}^{s} b_k \left( \Pi_{i+\frac{1}{2}}^k + \Pi_{i-\frac{1}{2}}^k \right) = \\ & E_i^{n+\frac{1}{2}} - \Delta v_i^T B A^{-1} \Delta U_i + \frac{1}{2} \lambda_i \sum_{k=1}^{s} b_k \left( \Pi_{i+\frac{1}{2}}^k + \Pi_{i-\frac{1}{2}}^k \right) \end{aligned} \tag{19}$$

where in the bottom row the following notation $\Delta U_i^T = \left( U_i^{(1)} - U_i^n, ..., U_i^{(s)} - U_i^n \right)$, $\Delta v_i^T = \left( v_i^{(1)} - v_i^n, ..., v_i^{(s)} - v_i^n \right)$ and $B = \text{diag}(b_1, ..., b_s)$ has been introduced.

---

[3] An entropy conservative scheme is one for which the numerical flux verifies $\Delta v_{i+\frac{1}{2}}^T \tilde{F}_{i+\frac{1}{2}} = \Delta \Theta_{i+\frac{1}{2}}$ and thus the entropy production $\Pi_{i+\frac{1}{2}} = 0$. In the scalar case, the entropy conservative flux is unique for each choice of entropy function and can be computed as $\tilde{F}_{i+\frac{1}{2}} = \Delta \Theta_{i+\frac{1}{2}} / \Delta v_{i+\frac{1}{2}}$. As explained in [12], entropy stable schemes can be constructed by coupling an entropy conservative scheme to a suitable diffusion operator.

Eq. (19) is the final form of the entropy production. It is identical to Eq. (14), but it is written in a form that facilitates subsequent analysis. This is clearly appreciated in the symmetric case (systems (1) with symmetric Jacobians $f_U$, which includes the scalar case), where a quadratic entropy function $\eta = \frac{1}{2} U^T U$ can be chosen. With this choice, $H(v) = (\eta_{UU})^{-1} = 1$, $v = \eta_U^T = U$ and $E_i^{n+\frac{1}{2}} = \frac{1}{2}(U_i^{n+1} - U_i^n)^2$, and thus the entropy evolution equation (19) can be written as

$$\eta_i^{n+1} - \eta_i^n + \lambda_i \sum_{k=1}^{s} b_k \left( \Phi_{i+\frac{1}{2}}^k - \Phi_{i-\frac{1}{2}}^k \right) = -\frac{1}{2} \Delta U^T Q \Delta U + \frac{1}{2} \lambda_i \sum_{k=1}^{s} b_k \left( \Pi_{i+\frac{1}{2}}^k + \Pi_{i-\frac{1}{2}}^k \right) \quad (20)$$

where the temporal entropy production is now a quadratic form with matrix

$$Q = BA^{-1} + A^{-T}B - A^{-T}bb^T A^{-1}.$$

## 2.2. Conditions for local entropy stability

For those systems that possess an entropy function, a desirable property of any numerical approximation is that the corresponding numerical entropy scheme correctly reproduces the physical entropy conservation law. This means, at the very least, that the scheme should dissipate entropy (using conventions such that the entropy function is essentially the negative of the physical entropy and thus decreases in physically allowed processes). This property can be guaranteed, usually under a suitable CFL restriction, if the spatial scheme is entropy stable [11] and the time integration scheme dissipates entropy or produces at most a controllable amount of unphysical entropy. CFL bounds of this sort can be straightforwardly –though somewhat tediously– derived as in [12, 1].

We have already seen that the Backward-Euler scheme is unconditionally entropy stable for any entropy function when coupled to entropy stable spatial schemes. In the general case, entropy stability must be analyzed by looking at the total entropy production of the scheme, which for a generic RK-$s$ scheme is given by Eq. (19). The fully discrete scheme is (locally) entropy stable provided that $S_i \leq 0$. The last term of the bottom row (the entropy production of the spatial scheme) is negative provided that $b_k \geq 0$ and that the numerical flux is entropy stable. This latter condition can be guaranteed by choosing the following flux [12]

$$F_{i+\frac{1}{2}} = \tilde{F}_{i+\frac{1}{2}} - \frac{1}{2} D_{i+\frac{1}{2}} \Delta_{i+\frac{1}{2}} v \quad (21)$$

where $\tilde{F}_{i+\frac{1}{2}}$ is an entropy conservative flux, which we write in viscosity form as $\tilde{F}_{i+1/2} = \frac{1}{2}(F(U_{i+1}) + F(U_i) - \tilde{Q}_{i+1/2} \Delta_{i+1/2} v)$, and $D_{i+\frac{1}{2}}$ are symmetric, positive-definite dissipation matrices. With that choice, the entropy production of the spatial scheme is semi-negative definite

$$\frac{1}{2} \sum_{k=1}^{s} b_k \left( \Pi_{i+\frac{1}{2}}^k + \Pi_{i-\frac{1}{2}}^k \right) = -\frac{1}{4} \sum_{k=1}^{s} b_k \left( (\Delta_{i+\frac{1}{2}} v^{(k)})^T D_{i+\frac{1}{2}}^{(k)} \Delta_{i+\frac{1}{2}} v^{(k)} + (\Delta_{i-\frac{1}{2}} v^{(k)})^T D_{i-\frac{1}{2}}^{(k)} \Delta_{i-\frac{1}{2}} v^{(k)} \right) \leq 0 \quad (22)$$

In the symmetric case with quadratic entropy, the temporal entropy dissipation can be written as $-\frac{1}{2}\Delta U^T Q \Delta U$ (Eq. (20)), so if the matrix $Q = BA^{-1} + A^{-T}B - A^{-T}bb^T A^{-1}$ is non-negative definite, the fully discrete scheme is entropy stable regardless of the CFL number. The conditions $b_k \geq 0$ and $Q$ non-negative definite imply that the method is algebraically stable[4] [2, 31, 32] which in turn entails *AN*-, *BN*-, *A*- and *B*-stability [32] [33]. Hence, algebraically stable RK schemes are also locally quadratic-entropy stable. This result follows trivially from the fact that, for the square entropy function, entropy stability is equivalent to stability of the solution in the $L_2$ norm, which is itself equivalent to algebraic stability [34].

Various well-known methods such as Gauss and Radau IIA, for example, are algebraically stable ( [2], Theorem 359C). The SDIRK methods that we will consider, on the other hand, are not algebraically stable and, in fact, they can be shown to yield positive temporal entropy production (see Fig. 3). Even the quadratic entropy stability of the algebraically stable temporal schemes need not hold for more general entropies, as is seen to be the case in Fig. 5.

Conditions for entropy stability for generic entropies can be derived as follows. We start by seeking conditions for entropy stability of the temporal scheme. We can bound $E_i^{n+\frac{1}{2}}$ in terms of the maximum condition number $K$ of the mean inverse Hessians

$$\bar{H}_i^{n+\frac{1}{2}} = \int_{-1/2}^{1/2} H(\tfrac{1}{2}(v_i^{n+1} + v_i^n) + \xi(v_i^{n+1} - v_i^n))d\xi$$

and

$$\bar{H}_i^{(k)} = \int_{-1/2}^{1/2} H(\tfrac{1}{2}(v_i^{(k)} + v_i^n) + \xi(v_i^{(k)} - v_i^n))d\xi, \ (k = 1,...,s)$$

(in such a way that $K^{-1}I_N \leq \bar{H} \leq K I_N$ ) as (see [12], eq. (7.16))

$$E_i^{n+\frac{1}{2}} \leq \frac{K^3}{2}\left\|\Delta^{n+\frac{1}{2}}U_i\right\|^2 = \frac{K^3}{2}\Delta U_i^T A^{-T}bb^T A^{-1}\Delta U_i \qquad (23)$$

Likewise, we have

$$-\Delta v_i^T BA^{-1}\Delta U_i = -\Delta U_i^T \hat{H}_i^{-1} BA^{-1}\Delta U_i = -\tfrac{1}{2}\Delta U_i^T (\hat{H}_i^{-1}BA^{-1} + A^{-T}B\hat{H}_i^{-1})\Delta U_i \qquad (24)$$

where $\hat{H}_i^{-1} = \mathrm{diag}((\bar{H}_i^{(1)})^{-1},...,(\bar{H}_i^{(s)})^{-1})$. Combining (23) and (24), it follows that the temporal scheme dissipates entropy provided that the matrix $K^3 A^{-T}bb^T A^{-1} - \hat{H}_i^{-1}BA^{-1} - A^{-T}B\hat{H}_i^{-1}$ is negative-semidefinite. Notice that in the quadratic case the Hessians reduce to the identity matrix (with condition number $K = 1$) and the above matrix is then identical to $-Q$.

---

[4] A Runge-Kutta method is algebraically stable if the matrices $B$ and $M = BA + A^T B - bb^T$ are both non-negative definite ( [32], p. 275). If $A$ is non-singular, then $M$ and $Q$ are congruent, $Q = A^{-T}MA^{-1}$, and thus $M$ is non-negative definite iff $Q$ is.

If the temporal scheme produces entropy, it is still possible to attain entropy stability by compensating the temporal production with the spatial dissipation under a CFL condition. This can be clearly seen in Fig. 4 for the SDIRK2 scheme ( [7] and section 3) with quadratic entropy. CFL bounds of this sort have been derived by Tadmor [12], which showed that the Forward-Euler scheme coupled to entropy stable spatial schemes is locally entropy stable with a sufficiently small CFL value. The analysis has been extended in [35] to prove the global entropy stability of the FE scheme coupled to higher-order (TECNO) spatial schemes under a CFL bound and in [1] to explicit, higher-order RK schemes. In the present case, a CFL condition for entropy stability can be derived in the general case as follows. We start by bounding the temporal production as

$$E_i^{n+\frac{1}{2}} \leq \frac{K^3}{2} \left\| \Delta^{n+\frac{1}{2}} U_i \right\|^2 = \frac{K^3}{2} (\lambda_i)^2 \left\| \sum_{k=1}^{s} b_k \left( F_{i+\frac{1}{2}}^{(k)} - F_{i-\frac{1}{2}}^{(k)} \right) \right\|^2$$

Since $b_k \geq 0$ we can bound the flux differences using Jensen's inequality, yielding

$$E_i^{n+\frac{1}{2}} \leq \frac{K^3}{2} (\lambda_i)^2 \sum_{k=1}^{s} b_k \left\| F_{i+\frac{1}{2}}^{(k)} - F_{i-\frac{1}{2}}^{(k)} \right\|^2$$

Writing the flux differences in the above equation as

$$F_{i+\frac{1}{2}}^{(k)} - F_{i-\frac{1}{2}}^{(k)} = \frac{1}{2}\left[ (B_{i+\frac{1}{2}}^{(k)} - \tilde{Q}_{i+\frac{1}{2}}^{(k)} - D_{i+\frac{1}{2}}^{(k)})\Delta_{i+1/2} v^{(k)} + (B_{i-\frac{1}{2}}^{(k)} + \tilde{Q}_{i-\frac{1}{2}}^{(k)} + D_{i-\frac{1}{2}}^{(k)})\Delta_{i-1/2} v^{(k)} \right]$$

($k = 1,\ldots,s$), where $\Delta_{i+1/2} v^{(k)} = v_{i+1}^{(k)} - v_i^{(k)}$ and the (symmetric) Jacobians $B_{i+\frac{1}{2}}^{(k)}$ are defined in the mean value sense as $F(U_{i+1}^{(k)}(v_{i+1}^{(k)})) - F(U_i^{(k)}(v_i^{(k)})) = B_{i+\frac{1}{2}}^{(k)} \cdot (v_{i+1}^{(k)} - v_i^{(k)})$; we get

$$E_i^{n+\frac{1}{2}} \leq \frac{3K^3}{4} (\lambda_i)^2 \sum_{k=1}^{s} b_k \left[ (\Delta_{i+\frac{1}{2}} v^{(k)})^T \left( (B_{i+\frac{1}{2}}^{(k)})^2 + (\tilde{Q}_{i+\frac{1}{2}}^{(k)})^2 + (D_{i+\frac{1}{2}}^{(k)})^2 \right) \Delta_{i+\frac{1}{2}} v^{(k)} \right.$$
$$\left. + (\Delta_{i-\frac{1}{2}} v^{(k)})^T \left( (B_{i-\frac{1}{2}}^{(k)})^2 + (\tilde{Q}_{i-\frac{1}{2}}^{(k)})^2 + (D_{i-\frac{1}{2}}^{(k)})^2 \right) \Delta_{i-\frac{1}{2}} v^{(k)} \right]$$

The remaining term can be bounded as

$$-\Delta v_i^T BA^{-1}\Delta U_i = \Delta v_i^T BR_i \leq \left|\Delta v_i^T BR_i\right| = \left|\sum_{k=1}^{s} b_k (U_i^{(k)} - U_i^n)^T \cdot (\bar{H}_i^{(k)})^{-1} \cdot R_i^{(k)}\right| =$$

$$\left|\sum_{k=1}^{s}\sum_{j=1}^{s} b_k A_{kj} (R_i^{(j)})^T \cdot (\bar{H}_i^{(k)})^{-1} \cdot R_i^{(k)}\right| \leq K \left|\sum_{k=1}^{s}\sum_{j=1}^{s} b_k A_{kj} (R_i^{(j)})^T \cdot R_i^{(k)}\right| \leq K \sum_{\alpha=1}^{N} \|BA\|\|R_i^\alpha\|^2 =$$

$$K(\lambda_i)^2 \sum_{\alpha=1}^{N}\sum_{k=1}^{s} \|BA\|\|(F_{i+\frac{1}{2}}^{(k)} - F_{i-\frac{1}{2}}^{(k)})^\alpha\|^2 = K(\lambda_i)^2 \sum_{k=1}^{s} \|BA\|\|F_{i+\frac{1}{2}}^{(k)} - F_{i-\frac{1}{2}}^{(k)}\|^2 \leq$$

$$K(\lambda_i)^2 \sum_{k=1}^{s} \|BA\|\frac{3}{2}\bigg[ (\Delta_{i+\frac{1}{2}}v^{(k)})^T \left((B_{i+\frac{1}{2}}^{(k)})^2 + (\tilde{Q}_{i+\frac{1}{2}}^{(k)})^2 + (D_{i+\frac{1}{2}}^{(k)})^2\right)\Delta_{i+\frac{1}{2}}v^{(k)}$$

$$+(\Delta_{i-\frac{1}{2}}v^{(k)})^T \left((B_{i-\frac{1}{2}}^{(k)})^2 + (\tilde{Q}_{i-\frac{1}{2}}^{(k)})^2 + (D_{i-\frac{1}{2}}^{(k)})^2\right)\Delta_{i-\frac{1}{2}}v^{(k)} \bigg]$$

In all, the temporal entropy production of the scheme obeys the following bound

$$E_i^{n+\frac{1}{2}} - \Delta v_i^T BA^{-1}\Delta U_i \leq$$
$$\frac{3}{2}K(\lambda_i)^2 \sum_{k=1}^{s}\left(\frac{K^2 b_k}{2} + \|BA\|\right)\bigg[ (\Delta_{i+\frac{1}{2}}v^{(k)})^T \left((B_{i+\frac{1}{2}}^{(k)})^2 + (\tilde{Q}_{i+\frac{1}{2}}^{(k)})^2 + (D_{i+\frac{1}{2}}^{(k)})^2\right)\Delta_{i+\frac{1}{2}}v^{(k)} \quad (25)$$
$$+(\Delta_{i-\frac{1}{2}}v^{(k)})^T \left((B_{i-\frac{1}{2}}^{(k)})^2 + (\tilde{Q}_{i-\frac{1}{2}}^{(k)})^2 + (D_{i-\frac{1}{2}}^{(k)})^2\right)\Delta_{i-\frac{1}{2}}v^{(k)} \bigg]$$

Using (25), the total entropy production (19) can be bounded as

$$S_i \leq \frac{\lambda_i}{4}\sum_{k=1}^{s}(\Delta_{i+\frac{1}{2}}v^{(k)})^T \Xi_{i+\frac{1}{2}}^{(k)}\Delta_{i+\frac{1}{2}}v^{(k)} + \frac{\lambda_i}{4}\sum_{k=1}^{s}(\Delta_{i-\frac{1}{2}}v^{(k)})^T \Xi_{i-\frac{1}{2}}^{(k)}\Delta_{i-\frac{1}{2}}v^{(k)} \quad (26)$$

where $\Xi_{i\pm\frac{1}{2}}^{(k)} = -b_k D_{i\pm\frac{1}{2}}^{(k)} + 6K\lambda_i\left(\frac{K^2 b_k}{2} + \|BA\|\right)\left((B_{i\pm\frac{1}{2}}^{(k)})^2 + (\tilde{Q}_{i\pm\frac{1}{2}}^{(k)})^2 + (D_{i\pm\frac{1}{2}}^{(k)})^2\right)$. It follows from (26) that the RK-s scheme is entropy stable, $S_i \leq 0$, provided that $\lambda_i$ is sufficiently small that $\Xi_{i\pm\frac{1}{2}}^{(k)} \leq 0$, or

$$b_k D_{i\pm\frac{1}{2}}^{(k)} \geq 6K\lambda_i\left(\frac{K^2 b_k}{2} + \|BA\|\right)\left((B_{i\pm\frac{1}{2}}^{(k)})^2 + (\tilde{Q}_{i\pm\frac{1}{2}}^{(k)})^2 + (D_{i\pm\frac{1}{2}}^{(k)})^2\right) \quad (27)$$

(we obviously need $b_k > 0$) for all $k$. The above bound is actually not too sharp, as among other things it assumes that $\Delta v_i^T BA^{-1}\Delta U_i$ is maximally negative, which need not be the case. We can confirm this idea with the Burgers equation with quadratic entropy and the SDIRK2 scheme. Now $K=1$, $v=U$, $B_{i+\frac{1}{2}} = \frac{1}{2}(U_{i+1} + U_i)$ and we pick the entropy stable flux

$$F_{i+\frac{1}{2}}^{ES} = \frac{1}{6}(U_{i+1}^2 + U_{i+1}U_i + U_i^2) - \mu(U_{i+1} - U_i) \quad (28)$$

(where $\mu > 0$) in such a way that $\tilde{Q}_{i+\frac{1}{2}} = \frac{1}{6}(U_{i+1} - U_i)$ and $D_{i+\frac{1}{2}} = 2\mu$. Hence

$$\min \lambda_i \leq \frac{\frac{2-\sqrt{2}}{2}\mu}{3\left(\frac{2-\sqrt{2}}{4}+\|BA\|\right)\max\left(\frac{1}{4}(U_{i+1}+U_i)^2+\frac{1}{36}(U_{i+1}-U_i)^2+4\mu^2\right)} \quad (29)$$

With the data of the example in section 3 we get, from (29)

$$\min \lambda_i \leq \frac{\frac{2-\sqrt{2}}{6}\mu}{\left(\frac{3-2\sqrt{2}}{2}+\|BA\|\right)\max\left(\frac{1}{4}(U_{i+1}+U_i)^2+\frac{1}{36}(U_{i+1}-U_i)^2+4\mu^2\right)} \approx 0.022$$

while numerical testing (Fig. 4) indicates a 50× larger CFL bound for entropy stability ($1.0 \leq \lambda_{\max} < 1.5$).

### 2.3. The linear symmetric case and stability

As we have just experienced, establishing general results about the sign of (14) or (19) is hard, but a simplification can be obtained by considering the symmetric case with entropy function $\eta = \frac{1}{2}U^T U$ for the linear semidiscrete problem:

$$U_t = LU \quad (30)$$

where the operator $L$ is constant in time. We will assume that (30) corresponds to the discretization of a scalar or symmetric problem such that it admits a quadratic entropy function. The evolution equation for the global quadratic entropy (or, equivalently, the energy) is

$$\frac{d}{dt}\left(\tfrac{1}{2}U^T U\right) = U^T U_t = U^T LU = \tfrac{1}{2}U^T(L-L^T)U + \tfrac{1}{2}U^T(L+L^T)U = \tfrac{1}{2}U^T(L+L^T)U \quad (31)$$

We can thus identify $U^T(L+L^T)U/2 = \langle U,(L+L^T)U\rangle = 2\langle U, LU\rangle$ with the net entropy production of the spatial scheme. Notice that if $L$ is skew-symmetric, then the scheme is (globally) entropy conservative, and if $L$ is semi-negative ($L+L^T \leq 0$) then the scheme is (globally) entropy stable. Coupling (30) to a time-advancement scheme results in a modified entropy balance equation. Taking into account that for quadratic entropies $v = U$ and $E_i^{n+\frac{1}{2}} = \tfrac{1}{2}(U_i^{n+1}-U_i^n)^T(U_i^{n+1}-U_i^n)$, the global entropy evolution equation (19) reads

$$\begin{aligned}
&\frac{1}{2}(U^{n+1})^T(U^{n+1}) - \frac{1}{2}(U^n)^T(U^n) = \frac{1}{2}\|U^{n+1}\|^2 - \frac{1}{2}\|U^n\|^2 = \\
&\frac{1}{2}\|U^{n+1}-U^n\|^2 - \Delta U^T BA^{-1}\Delta U + \frac{1}{2}\sum_{k=1}^s b_k \Delta t \langle U^{(k)},(L+L^T)U^{(k)}\rangle = \\
&\frac{1}{2}\Delta U^T(A^{-T}bb^T A^{-1} - BA^{-1} - A^{-T}B)\Delta U + \frac{1}{2}\sum_{k=1}^s b_k \Delta t \langle U^{(k)},(L+L^T)U^{(k)}\rangle = \\
&-\frac{1}{2}\Delta U^T Q \Delta U + \frac{1}{2}\sum_{k=1}^s b_k \Delta t \langle U^{(k)},(L+L^T)U^{(k)}\rangle
\end{aligned} \quad (32)$$

It is important to stress that (32) is actually valid for the energy variation of *any* linear system. Only in the symmetric case it also represents the entropy evolution of the scheme.

Notice that if $b_k \geq 0$ and $Q = BA^{-1} + A^{-T}B - A^{-T}bb^T A^{-1}$ is non-negative definite (and thus the method is algebraically stable) then the method is strongly stable ($\|U^{n+1}\|^2 - \|U^n\|^2 \leq 0$) for general semi-negative operators $L$ obeying $L+L^T \leq 0$.

### 2.4. Examples

In order to check the above formulas and to pave the way for subsequent numerical testing, we now review several examples of RK schemes and compute their corresponding entropy evolution equations.

(a) The diagonally implicit (DIRK) $s$-stage scheme with tableau

$$\begin{array}{c|cccc} b_1 & b_1 & & & \\ b_1+b_2 & b_1 & b_2 & & \\ \cdots & \cdots & \cdots & & \\ 1 & b_1 & b_2 & \cdots & b_s \\ \hline & b_1 & b_2 & \cdots & b_s \end{array}$$

(and $b_1 b_2 \cdots b_s \neq 0$) corresponds to $s$ consecutive (and independent) Backward-Euler steps. Since

$$A^{-1} = \begin{pmatrix} \frac{1}{b_1} & 0 & 0 & \cdots & 0 \\ -\frac{1}{b_1} & \frac{1}{b_2} & 0 & \cdots & 0 \\ 0 & -\frac{1}{b_2} & \frac{1}{b_3} & \cdots & 0 \\ \vdots & \vdots & \ddots & \ddots & \vdots \\ 0 & 0 & \cdots & -\frac{1}{b_{s-1}} & \frac{1}{b_s} \end{pmatrix}$$

the entropy production for this scheme is, from (14),

$$\eta_i^{n+1} - \eta_i^n + \lambda_i \sum_{k=1}^{s} b_k \left( \Phi_{i+\frac{1}{2}}^k - \Phi_{i-\frac{1}{2}}^k \right) =$$

$$-B_i^{n+\frac{1}{2}} + \sum_{k=1}^{s}\sum_{j=1}^{s} b_k A_{kj}^{-1}(v_i^{n+1} - v_i^{(k)})^T (U_i^{(j)} - U_i^n) + \frac{1}{2}\lambda_i \sum_{k=1}^{s} b_k \left( \Pi_{i+\frac{1}{2}}^k + \Pi_{i-\frac{1}{2}}^k \right) =$$

$$-\sum_{k=1}^{s} B_i^{k-1,k} + \frac{1}{2}\lambda_i \sum_{k=1}^{s} b_k \left( \Pi_{i+\frac{1}{2}}^k + \Pi_{i-\frac{1}{2}}^k \right)$$

(where $-B_i^{k-1,k} = \eta_i^{(k)} - \eta_i^{(k-1)} - (v_i^{(k)})^T (U_i^{(k)} - U_i^{(k-1)})$), which corresponds simply to $s$ copies of the single-step Backward-Euler entropy production term (16), as expected.

(b) Crank-Nicolson (Trapezoidal scheme)

$$U_i^{n+1} = U_i^n - \frac{1}{2} R_i^n - \frac{1}{2} R_i^{n+1} = U_i^n - \bar{R}_i^{n+\frac{1}{2}}$$

This is a 2-stage Runge-Kutta scheme with tableau

|   | 0   | 0   |
|---|-----|-----|
| 1 | 1/2 | 1/2 |
|   | 1/2 | 1/2 |

Notice that since $A$ is not invertible and $\tilde{A}$ has rank 1, eq. (14) or (19) cannot be used here. An entropy evolution equation can be still obtained by combining (10) and (17):

$$\eta_i^{n+1} - \eta_i^n = \tfrac{1}{2}(E_i^{n+\frac{1}{2}} - B_i^{n+\frac{1}{2}}) + (\overline{v}_i^{n+\frac{1}{2}})^T \Delta U_i^{n+\frac{1}{2}} = \tfrac{1}{2}(E_i^{n+\frac{1}{2}} - B_i^{n+\frac{1}{2}}) - (\overline{v}_i^{n+\frac{1}{2}})^T \overline{R}_i^{n+\frac{1}{2}}$$

(1) quadratic entropy: $\overline{v}_i^{n+\frac{1}{2}} = v(\overline{U}_i^{n+\frac{1}{2}}) = \overline{U}_i^{n+\frac{1}{2}} \equiv v_i^{n+\frac{1}{2}}$ and $E_i^{n+\frac{1}{2}} - B_i^{n+\frac{1}{2}} = 0$. Hence,

$$\eta_i^{n+1} - \eta_i^n = -(v_i^{n+\frac{1}{2}})^T \overline{R}_i^{n+\frac{1}{2}}$$

(2) quadratic entropy + linear fluxes: additionally to the above, $\overline{R}_i^{n+\frac{1}{2}} = R_i(\overline{U}^{n+\frac{1}{2}})$ and thus

$$\eta_i^{n+1} - \eta_i^n = -(v_i^{n+\frac{1}{2}})^T R_i^{n+\frac{1}{2}}(\overline{U}_i^{n+\frac{1}{2}}) = -\lambda_i (\Delta \Phi_i^{n+\frac{1}{2}} - \Pi_i^{n+\frac{1}{2}})$$

(the temporal scheme is entropy conservative in this case).

(3) general case

$$\eta_i^{n+1} - \eta_i^n = \tfrac{1}{2}(E_i^{n+\frac{1}{2}} - B_i^{n+\frac{1}{2}}) - (\overline{v}_i^{n+\frac{1}{2}})^T \overline{R}_i^{n+\frac{1}{2}} =$$

$$\tfrac{1}{2}(E_i^{n+\frac{1}{2}} - B_i^{n+\frac{1}{2}}) - \tfrac{1}{2}((v_i^{n+1})^T R_i^{n+1} + (v_i^n)^T R_i^n) + \tfrac{1}{4}(v_i^{n+1} - v_i^n)^T (R_i^{n+1} - R_i^n) =$$

$$-\tfrac{1}{2}\lambda(\Phi_{i+\frac{1}{2}}^{n+1} - \Phi_{i-\frac{1}{2}}^{n+1} + \Phi_{i+\frac{1}{2}}^n - \Phi_{i-\frac{1}{2}}^n) + \tfrac{1}{4}\lambda(\Pi_{i+\frac{1}{2}}^{n+1} + \Pi_{i-\frac{1}{2}}^{n+1} + \Pi_{i+\frac{1}{2}}^n + \Pi_{i-\frac{1}{2}}^n) + \tfrac{1}{2}(E_i^{n+\frac{1}{2}} - B_i^{n+\frac{1}{2}})$$

$$+ \tfrac{1}{4}(v_i^{n+1} - v_i^n)^T (R_i^{n+1} - R_i^n)$$

or, rearranging,

$$\eta_i^{n+1} - \eta_i^n + \tfrac{1}{2}\lambda(\Phi_{i+\frac{1}{2}}^{n+1} - \Phi_{i-\frac{1}{2}}^{n+1} + \Phi_{i+\frac{1}{2}}^n - \Phi_{i-\frac{1}{2}}^n) =$$

$$\underbrace{\tfrac{1}{4}\lambda(\Pi_{i+\frac{1}{2}}^{n+1} + \Pi_{i-\frac{1}{2}}^{n+1} + \Pi_{i+\frac{1}{2}}^n + \Pi_{i-\frac{1}{2}}^n)}_{S_i^{(x)}} + \underbrace{\tfrac{1}{2}(E_i^{n+\frac{1}{2}} - B_i^{n+\frac{1}{2}}) + \tfrac{1}{4}\Delta t(v_i^{n+1} - v_i^n)^T (R_i^{n+1} - R_i^n)}_{S_i^{(t)}}$$

(c) Generalized Crank-Nicolson (Tadmor [12])

$$U_i^{n+1} = U_i^n - R_i(U(\tilde{v}^{n+\frac{1}{2}}))$$

where $\tilde{v}^{n+\frac{1}{2}} = \int_{-1/2}^{1/2} v(\overline{U}^{n+\frac{1}{2}} + \Delta U^{n+\frac{1}{2}}\xi)d\xi$. This is not a RK scheme, but it is included here for completeness. Now, since

$$\eta(U^{n+1}) - \eta(U^n) =$$

$$\int_{U^n}^{U^{n+1}} \eta_U dU = \int_{U^n}^{U^{n+1}} v^T(U)dU = \left[\int_{-1/2}^{1/2} v^T(\overline{U}^{n+\frac{1}{2}} + \Delta U^{n+\frac{1}{2}}\xi)d\xi\right]\Delta U^{n+\frac{1}{2}} = (\tilde{v}^{n+\frac{1}{2}})^T \Delta U^{n+\frac{1}{2}}$$

it turns out that the discrete entropy evolution equation is simply

$$\eta_i^{n+1} - \eta_i^n = (\tilde{v}_i^{n+\frac{1}{2}})^T \Delta U_i^{n+\frac{1}{2}} = -\Delta t (\tilde{v}_i^{n+\frac{1}{2}})^T R_i (U(\tilde{v}^{n+\frac{1}{2}})) = -\lambda \left( \Delta \Phi_i(\tilde{v}^{n+\frac{1}{2}}) - \Pi_i(\tilde{v}^{n+\frac{1}{2}}) \right)$$

and thus the temporal scheme does not produce entropy. This scheme is identical to Crank-Nicolson in the symmetric case (with quadratic entropy) and linear fluxes (cf. case (b.3) above).

(d) 4$^{th}$ order, 2-stage Gauss scheme [2] with Tableau

$$\begin{array}{c|cc} \frac{1}{2} - \frac{\sqrt{3}}{6} & \frac{1}{4} & \frac{1}{4} - \frac{\sqrt{3}}{6} \\ \frac{1}{2} + \frac{\sqrt{3}}{6} & \frac{1}{4} + \frac{\sqrt{3}}{6} & \frac{1}{4} \\ \hline & \frac{1}{2} & \frac{1}{2} \end{array}$$

Since

$$A^{-1} = \begin{pmatrix} 3 & -3+2\sqrt{3} \\ -3-2\sqrt{3} & 3 \end{pmatrix}$$

the entropy evolution equation is

$$\eta_i^{n+1} - \eta_i^n + \frac{1}{2} \lambda_i \sum_{k=1}^{2} \left( \Phi_{i+\frac{1}{2}}^k - \Phi_{i-\frac{1}{2}}^k \right) =$$

$$E_i^{n+\frac{1}{2}} - \frac{1}{2}((v_i^{(1)} - v_i^n)^T \; (v_i^{(2)} - v_i^n)^T) \begin{pmatrix} 3 & -3+2\sqrt{3} \\ -3-2\sqrt{3} & 3 \end{pmatrix} \begin{pmatrix} U_i^{(1)} - U_i^n \\ U_i^{(2)} - U_i^n \end{pmatrix} + \frac{\lambda_i}{4} \sum_{k=1}^{2} \left( \Pi_{i+\frac{1}{2}}^k + \Pi_{i-\frac{1}{2}}^k \right)$$

It can be checked that $Q = 0$ for this method [32]. Hence, in the symmetric case with quadratic entropy the above evolution equation reduces to

$$\frac{1}{2}(U_i^{n+1})^T U_i^{n+1} - \frac{1}{2}(U_i^n)^T U_i^n + \frac{1}{2} \lambda_i \sum_{k=1}^{2} \left( \Phi_{i+\frac{1}{2}}^k - \Phi_{i-\frac{1}{2}}^k \right) = \frac{1}{4} \lambda_i \sum_{k=1}^{2} \left( \Pi_{i+\frac{1}{2}}^k + \Pi_{i-\frac{1}{2}}^k \right)$$

Thus, the $s = 2$ Gauss temporal scheme does not produce entropy in this case and, when coupled to entropy stable spatial operators, the resulting fully discrete scheme is quadratic entropy stable (and also strongly stable for general linear semi-negative operators). Furthermore, these results hold for all Gauss methods, which have $Q = 0$ [32].

(e) 3$^{th}$ order, 2-stage Radau IIA scheme [2] with Tableau

$$\begin{array}{c|cc} \frac{1}{3} & \frac{5}{12} & -\frac{1}{12} \\ 1 & \frac{3}{4} & \frac{1}{4} \\ \hline & \frac{3}{4} & \frac{1}{4} \end{array}$$

Here

$$A^{-1} = \begin{pmatrix} \frac{3}{2} & \frac{1}{2} \\ -\frac{9}{2} & \frac{5}{2} \end{pmatrix}$$

and, thus,

$$\eta_i^{n+1} - \eta_i^n + \lambda_i \sum_{k=1}^{2} b_k \left( \Phi_{i+\frac{1}{2}}^k - \Phi_{i-\frac{1}{2}}^k \right) =$$

$$E_i^{n+\frac{1}{2}} - ((v_i^{(1)} - v_i^n)^T \ (v_i^{(2)} - v_i^n)^T) \begin{pmatrix} \frac{9}{8} & \frac{3}{8} \\ -\frac{9}{8} & \frac{5}{8} \end{pmatrix} \begin{pmatrix} U_i^{(1)} - U_i^n \\ U_i^{(2)} - U_i^n \end{pmatrix} + \frac{1}{2} \lambda_i \sum_{k=1}^{2} b_k \left( \Pi_{i+\frac{1}{2}}^k + \Pi_{i-\frac{1}{2}}^k \right)$$

The matrix $Q = \begin{pmatrix} \frac{9}{4} & \frac{-3}{4} \\ \frac{-3}{4} & \frac{1}{4} \end{pmatrix}$ for this method, so the entropy evolution equation in the symmetric case is

$$\frac{1}{2}(U_i^{n+1})^T U_i^{n+1} - \frac{1}{2}(U_i^n)^T U_i^n + \lambda_i \sum_{k=1}^{2} b_k \left( \Phi_{i+\frac{1}{2}}^k - \Phi_{i-\frac{1}{2}}^k \right) =$$

$$-\frac{1}{2}((U_i^{(1)} - U_i^n)^T \ (U_i^{(2)} - U_i^n)^T) \begin{pmatrix} \frac{9}{4} & \frac{-3}{4} \\ \frac{-3}{4} & \frac{1}{4} \end{pmatrix} \begin{pmatrix} U_i^{(1)} - U_i^n \\ U_i^{(2)} - U_i^n \end{pmatrix} + \frac{1}{2} \lambda_i \sum_{k=1}^{2} b_k \left( \Pi_{i+\frac{1}{2}}^k + \Pi_{i-\frac{1}{2}}^k \right)$$

Since $Q$ is positive-semidefinite, the $s = 2$ Radau IIA method coupled to entropy stable spatial schemes is locally entropy stable and also strongly stable for general linear semi-negative operators. Furthermore, the above (entropy) stability remarks hold for all Radau IIA methods, which are algebraically stable (and thus have positive-semidefinite $Q$).

## 3. Numerical experiments

In this section, we will use the above formulas to compute the entropy production in two simple numerical test cases. Temporal integration is performed with the Backward-Euler and two-stage Gauss and Radau IIA schemes described in Section 2, as well as with the 6[th]-order, 3-stage, $A$-stable Gauss scheme with tableau

$$\begin{array}{c|ccc} \frac{1}{2} - \frac{\sqrt{15}}{10} & \frac{5}{36} & \frac{2}{9} - \frac{\sqrt{15}}{15} & \frac{5}{36} - \frac{\sqrt{15}}{30} \\ \frac{1}{2} & \frac{5}{36} + \frac{\sqrt{15}}{30} & \frac{2}{9} & \frac{5}{36} - \frac{\sqrt{15}}{24} \\ \frac{1}{2} + \frac{\sqrt{15}}{10} & \frac{5}{36} + \frac{\sqrt{15}}{30} & \frac{2}{9} + \frac{\sqrt{15}}{15} & \frac{5}{36} \\ \hline & \frac{5}{18} & \frac{8}{18} & \frac{5}{18} \end{array}$$

the 5[th]-order, 3-stage, $L$-stable Radau IIA scheme with tableau

$$\begin{array}{c|ccc} \frac{4-\sqrt{6}}{10} & \frac{88-7\sqrt{6}}{360} & \frac{296-169\sqrt{6}}{1800} & \frac{-2+3\sqrt{6}}{225} \\ \frac{4+\sqrt{6}}{10} & \frac{296+169\sqrt{6}}{1800} & \frac{88+7\sqrt{6}}{360} & \frac{-2-3\sqrt{6}}{225} \\ 1 & \frac{16-\sqrt{6}}{36} & \frac{16+\sqrt{6}}{36} & \frac{1}{9} \\ \hline & \frac{16-\sqrt{6}}{36} & \frac{16+\sqrt{6}}{36} & \frac{1}{9} \end{array}$$

and the 2[nd] and 3[rd]-order 2-stage and 3-stage $L$-stable SDIRK schemes [7] with tableaux

$$\begin{array}{c|cc}
1-\frac{1}{2}\sqrt{2} & 1-\frac{1}{2}\sqrt{2} & 0 \\
1 & \frac{1}{2}\sqrt{2} & 1-\frac{1}{2}\sqrt{2} \\
\hline
 & \frac{1}{2}\sqrt{2} & 1-\frac{1}{2}\sqrt{2}
\end{array}$$

(SDIRK2) and

$$\begin{array}{c|ccc}
\lambda & \lambda & 0 & 0 \\
\frac{1}{2}(1+\lambda) & \frac{1}{2}(1-\lambda) & \lambda & 0 \\
1 & \frac{1}{4}(-6\lambda^2+16\lambda-1) & \frac{1}{4}(6\lambda^2-20\lambda+5) & \lambda \\
\hline
 & \frac{1}{4}(-6\lambda^2+16\lambda-1) & \frac{1}{4}(6\lambda^2-20\lambda+5) & \lambda
\end{array}$$

(SDIRK3), where $\lambda = 0.4358665215$.

As our first test case, we consider the advection equation

$$u_t = u_x \tag{33}$$

in $-1 \leq x \leq 1$ with periodic boundary conditions and initial condition:

$$u(x,0) = \sin(\pi x) \tag{34}$$

With these choices, any (smooth) entropy function $\eta$ is globally constant, $\frac{d}{dt}\int_{-1}^{1} \eta dx = 0$. We focus on the quadratic entropy $\eta = \frac{1}{2}u^2$ for which $\int_{-1}^{1} \eta(t,x)dx = \frac{1}{2}$. A simple, entropy stable, spatial discretization is obtained by coupling the entropy conservative scheme $\tilde{F}_{i+\frac{1}{2}} = \Delta\Theta_{i+\frac{1}{2}}/\Delta v_{i+\frac{1}{2}}$ to the first-order dissipation operator $d_{i+\frac{1}{2}} = -\mu(v_{i+1} - v_i)$, with $\mu > 0$ [36]. Accordingly, eq. (33) is discretized on a uniform mesh with the following quadratic-entropy-stable scheme

$$u_t(x_j, t) = \frac{f_{j+\frac{1}{2}} - f_{j-\frac{1}{2}}}{\Delta x} \tag{35}$$

where

$$f_{j+\frac{1}{2}} = \frac{1}{2}(u_{j+1} + u_j) + \mu(u_{j+1} - u_j), \quad \mu \geq 0 \tag{36}$$

which reduces to the simple first-order upwind difference $du_j/dt = (u_{j+1} - u_j)/\Delta x$ when $\mu = 1/2$.

Fig. 1 shows a convergence study of the global error in the quadratic entropy $\left|\int_{-1}^{1} \eta(t,x)dx - \frac{1}{2}\right|$ at $t = 2$ for several RK schemes for the linear advection problem (33)-(36) with $\mu = 0$ and $n = 400$ cells in space. The results show that (1) the methods achieve their design order and (2) for a fixed time-step size, the errors decrease with the order $p$ of the scheme.

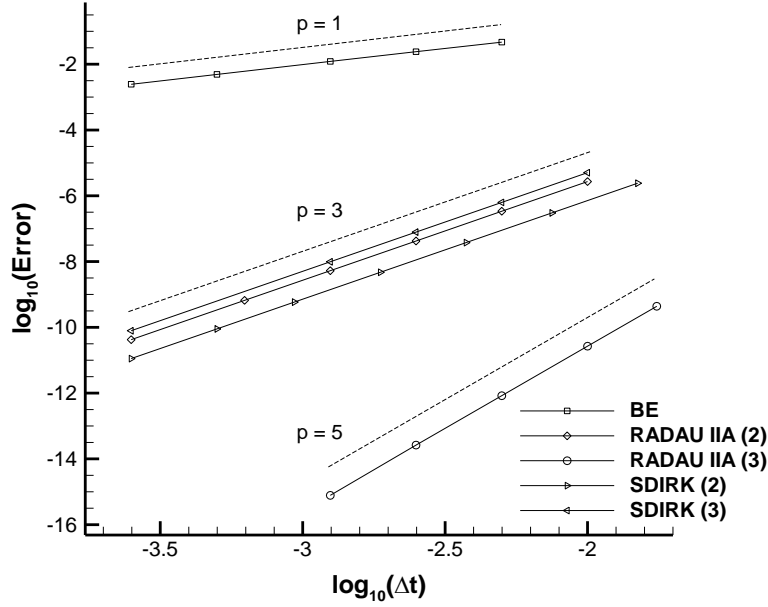

**Fig. 1** Linear advection at $t = 2$ with quadratic-entropy conservative scheme. Convergence of error in global entropy for various RK schemes.

For the next example we consider the inviscid Burgers equation

$$\frac{\partial U}{\partial t} + \frac{1}{2}\frac{\partial U^2}{\partial x} = 0$$

with initial data

$$U(x,0) = \begin{cases} 1.5 & (x \leq 0) \\ 0.5 & (x > 0) \end{cases}$$

which results in a right-moving shock with speed $s = 1$. We will compute the numerical solution of the above test case with the quadratic-entropy stable spatial scheme (28) as well as with the scheme

$$F^{ES}_{i+\frac{1}{2}} = \frac{1}{2}U_{i+1}U_i - \mu \frac{U_{i+1} - U_i}{U_{i+1}U_i} \tag{37}$$

(which is entropy stable for the logarithmic entropy $\eta = -\log(U)$ if $\mu > 0$).

Fig. 2 shows the total local quadratic and logarithmic entropy production $S_i$ at $t = 3$. We see that with the chosen value of the CFL number $\lambda = \Delta t / \Delta x$ all schemes are entropy stable. The separate contribution of the temporal scheme to the entropy production is examined in Fig. 3 (for the quadratic entropy) and Fig. 5 (for the logarithmic entropy). In the quadratic case, all schemes dissipate entropy (the Gauss schemes even being entropy-conservative) excepting the SDIRK schemes, which are seen to actually produce significant amounts of entropy at some locations, which must be compensated by the entropy dissipation of the spatial scheme in order to attain entropy stability. This is confirmed in Fig. 4, which plots the total (quadratic) entropy production of the SDIRK2 scheme for different values of the CFL number $\lambda = \Delta t / \Delta x$, showing that beyond a certain value of $\lambda$ the fully discrete scheme becomes entropy unstable.

Finally, Fig. 5 shows that in the logarithmic case, all schemes (excepting the BE scheme, which we already know to be unconditionally entropy stable), produce entropy. We also note from Fig. 3 and Fig. 5 that the order of each method is directly correlated with the size of the temporal entropy production: for a method of order $p$, the maximum value of the temporal entropy production is roughly of size $\sim 10^{-2p}$. On the other hand, the integrated temporal entropy production $\int S^{(t)} dx$ is shown to converge as $\Delta t^{p+1}$ in Fig. 6.

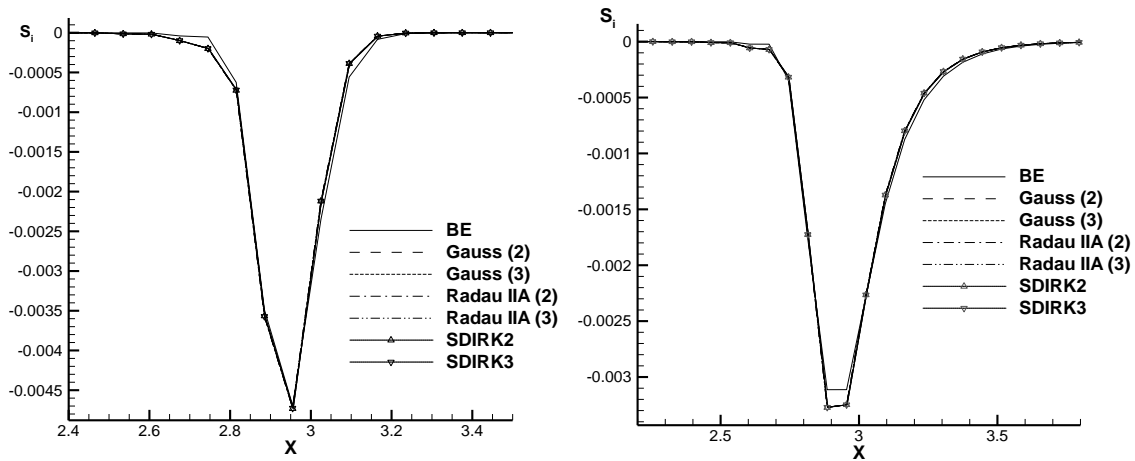

**Fig. 2 Burgers equation. Moving shock at $t = 3$. Total (local) entropy production for the quadratic (left) and logarithmic (right) entropy stable schemes.**

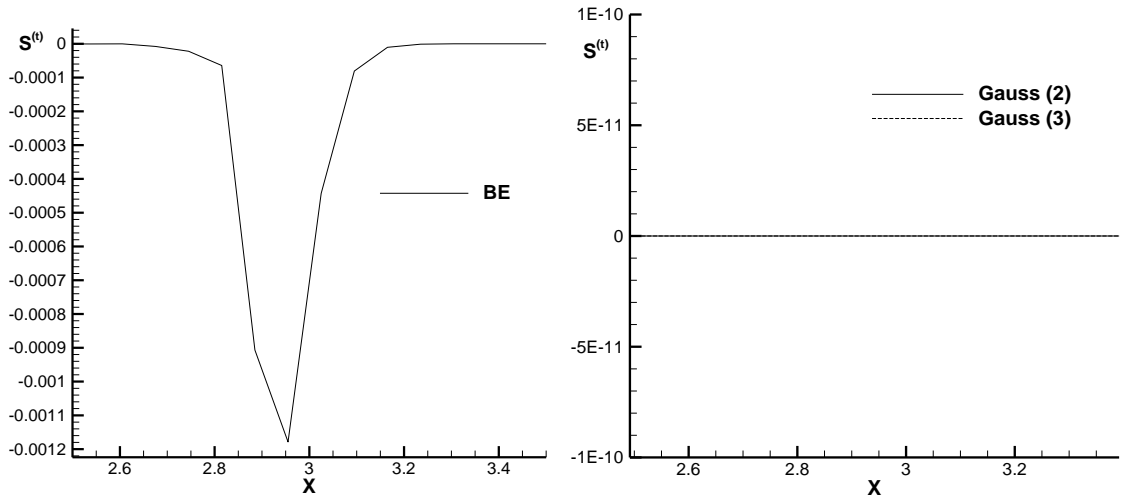

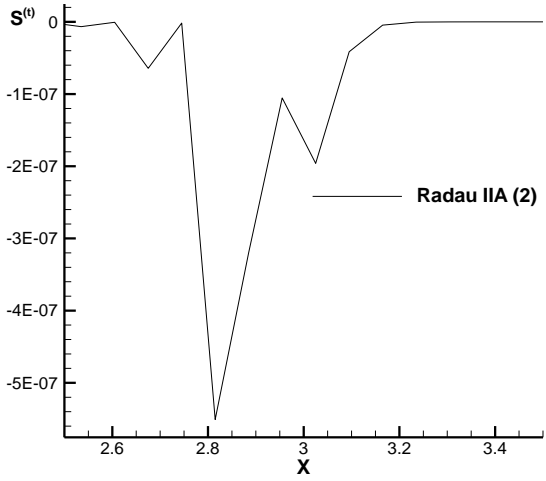 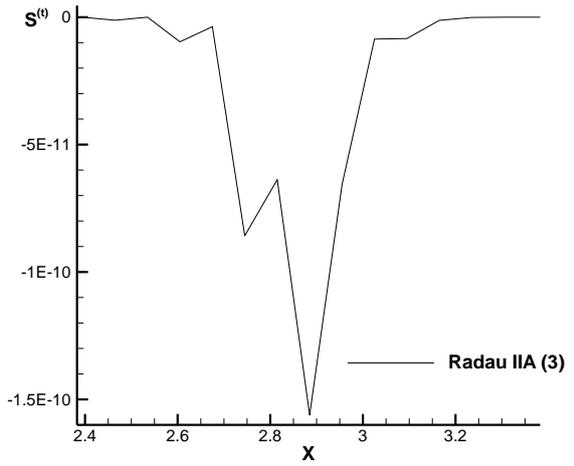
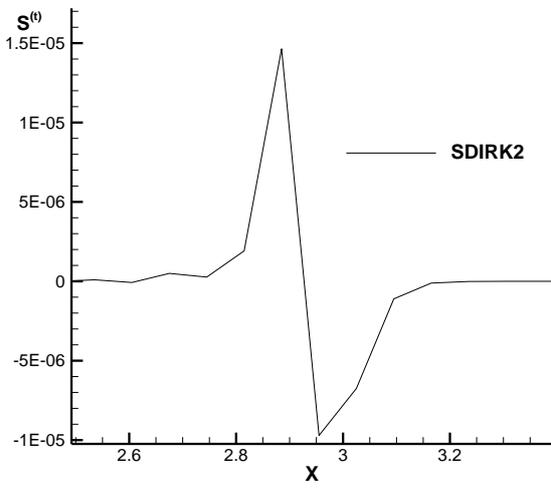 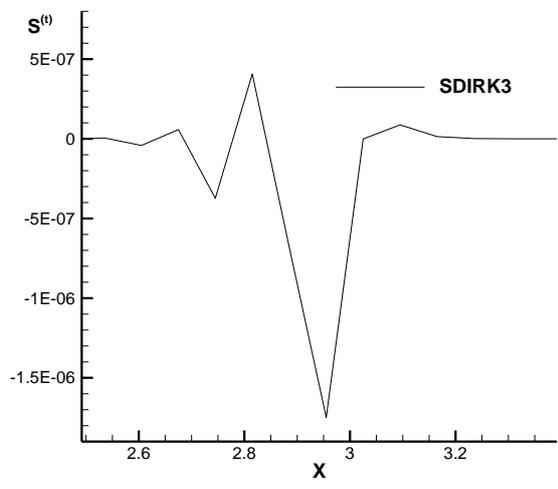

**Fig. 3 Burgers equation. Moving shock at $t = 3$. Local temporal entropy production for the quadratic entropy stable scheme.**

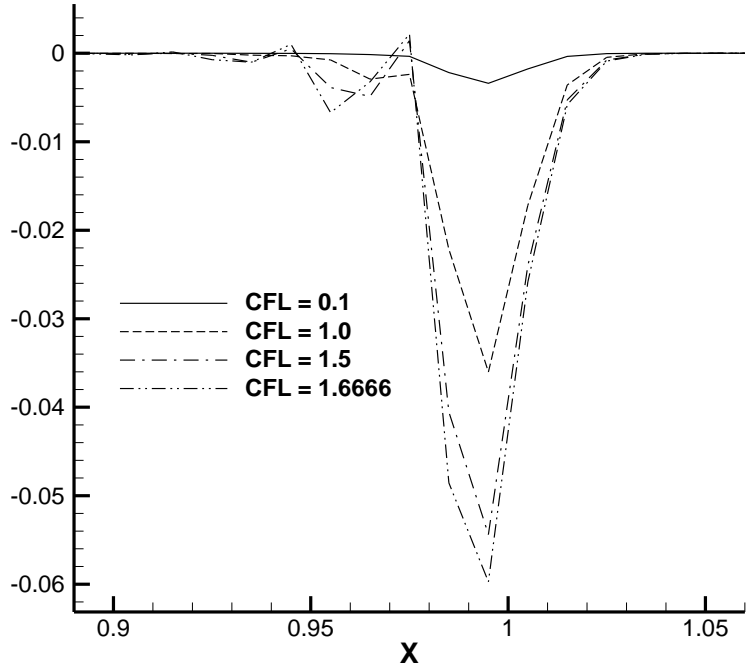

**Fig. 4** Burgers equation. Moving shock at $t = 3$. Total entropy production for the SDIRK2 scheme for the quadratic entropy stable scheme with different time-steps.

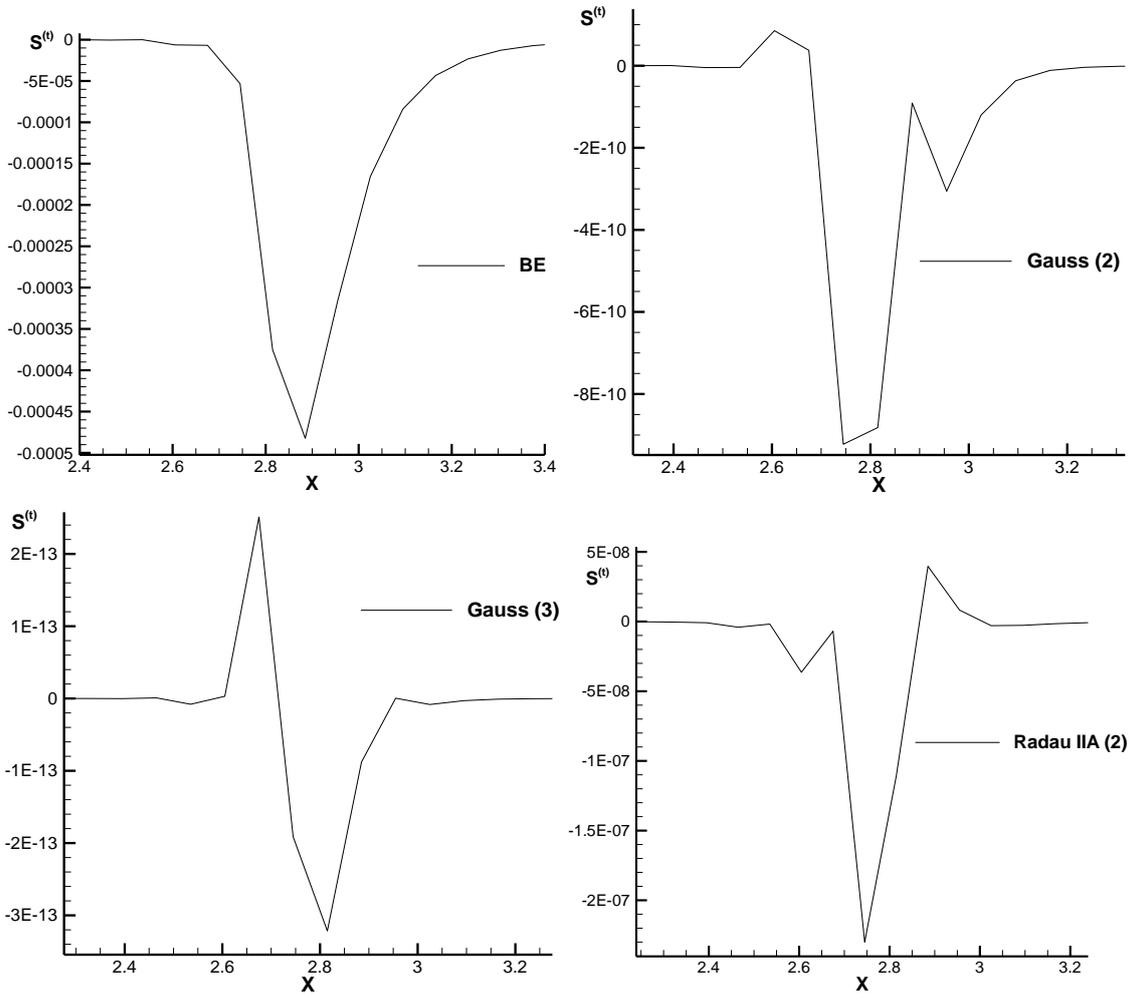

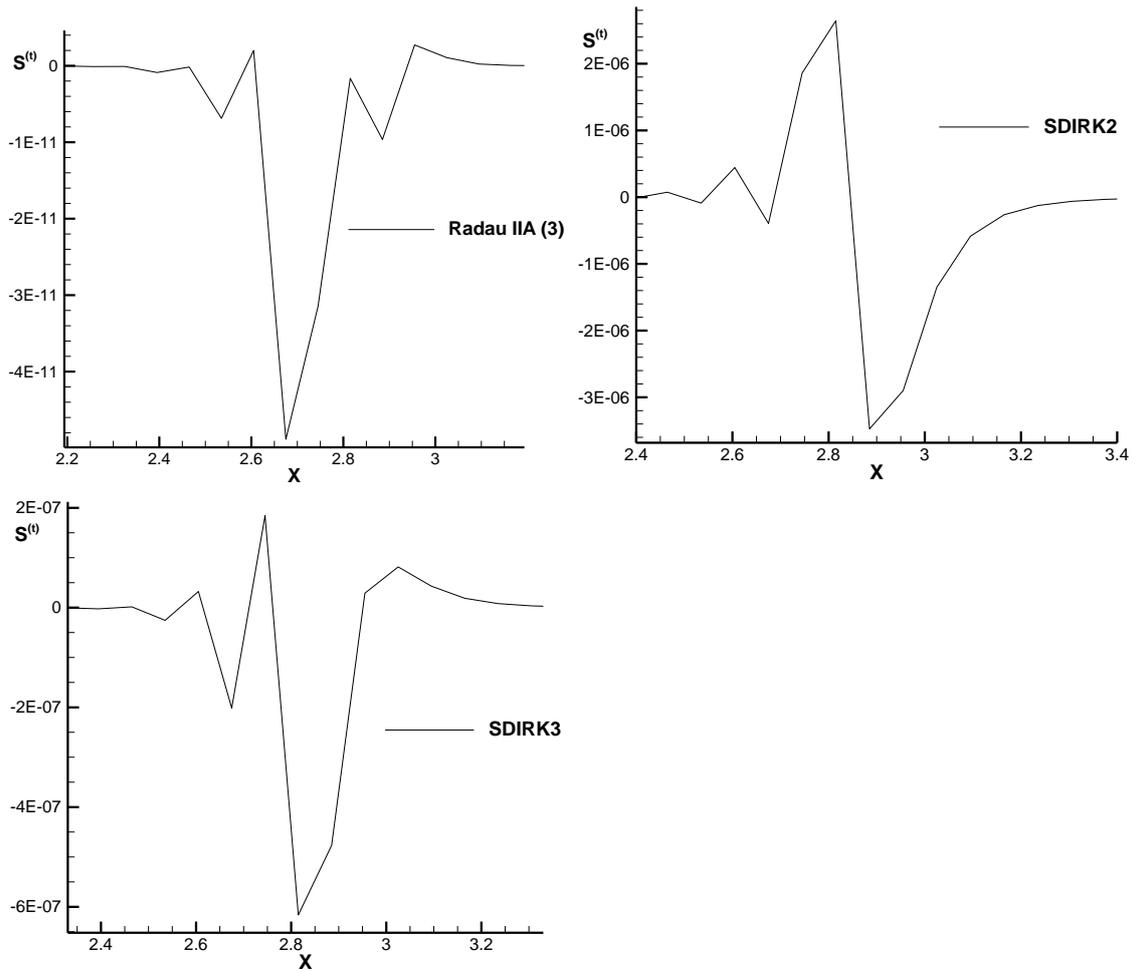

**Fig. 5 Burgers equation. Moving shock at *t* = 3. Local temporal entropy production at *t* =3 for the logarithmic entropy stable scheme.**

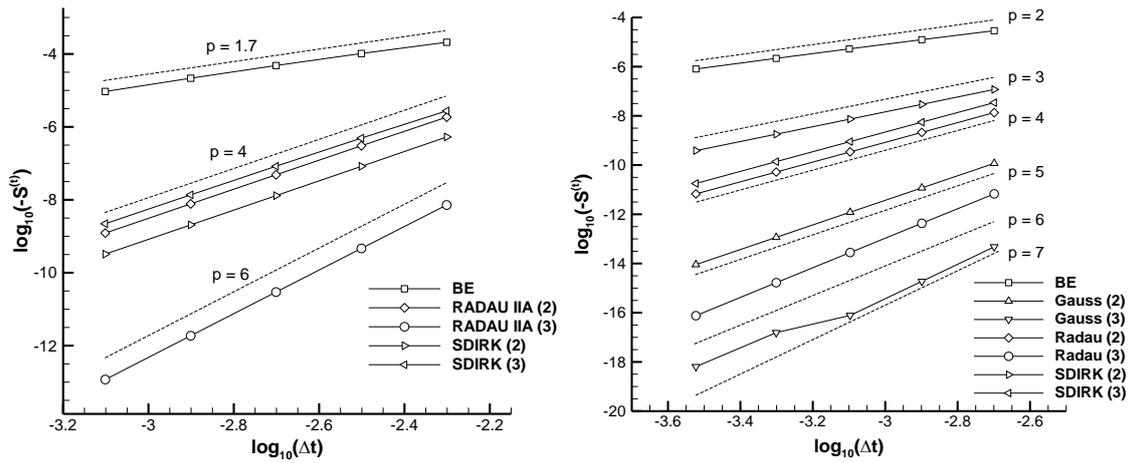

**Fig. 6 Burgers equation. Moving shock at *t* = 3. Convergence study of the integrated temporal entropy production for various RK schemes. Left: quadratic entropy. Right: logarithmic entropy.**

## 4. Conclusions

In this paper, we have derived the exact entropy evolution equation for fully discrete schemes arising from the finite-volume discretization of systems of conservation laws coupled to generic Runge-Kutta time integrators. The main results are summarized in Eqs. (13), (14) and (19). Using these representations, the conditions under which the fully discrete scheme is (locally) entropy stable have been examined, resulting in the general case in a CFL bound (27) which is however not particularly sharp.

In the particular case of the quadratic entropy function, which is relevant in the scalar case or for systems of conservation laws with symmetric flux Jacobians, algebraically stable RK schemes are stable in the $L_2$ norm and are thus locally temporally entropy stable. This latter result does not hold for more general entropies (as clearly shown in numerical testing), and it turns out that in the general case even fully implicit RK schemes produce entropy (two notable exceptions being the backward-Euler and the modified Crank-Nicolson schemes, which are unconditionally entropy stable –resp. entropy conservative– for any entropy.)

Finally, several, well-known implicit RK schemes have been examined. Gauss and Radau IIA methods are unconditionally quadratic-entropy stable, the former being actually quadratic entropy-conservative, while the examined SDIRK schemes produce entropy even in the quadratic case. All those assertions have been verified with numerical testing in both linear (advection of a sine wave) and non-linear (inviscid Burgers equation) cases.

## Acknowledgments

This work has been supported by the Spanish Ministry of Defence and INTA under the research program "Termofluidodinámica" (IGB99001).